\documentclass[a4paper,11pt]{article}
\usepackage{macros}
\usepackage{vmargin}

\setmarginsrb{2cm}{1.5cm}{2cm}{1cm}{0mm}{1cm}{0mm}{0.5cm}

\long\def\symbolfootnote[#1]#2{\begingroup%
\def\thefootnote{\fnsymbol{footnote}}\footnote[#1]{#2}\endgroup}

\begin{document}
\pagestyle{empty}

\begin{center}{\LARGE\bf Discussion of the Gear-Gupta-Leimkuhler method for impacting mechanical systems\symbolfootnote[2]{This is a preprint of a paper submitted to Multibody System Dynamics.}}\end{center} 

\begin{center}
\parbox{5in}{ \centering \textbf{Svenja Schoeder, Heinz Ulbrich, Thorsten Schindler}\\
  Institute of Applied Mechanics\\
	Technische Universit\"at M\"unchen\\
	Boltzmannstra\ss{}e 15\\
	85748 Garching, Germany\\
  {\tt\small thorsten.schindler@mytum.de}
}
\end{center}

\section*{Abstract}
In multibody simulation, the Gear-Gupta-Leimkuhler method for only persistent contacts enforces constraints on position and velocity level at the same time. It yields a robust numerical discretization of differential algebraic equations avoiding the drift-off effect. In this work, we carry over these benefits to impacting mechanical systems with unilateral constraints. For this kind of a mechanical system, adding the position level constraint to a timestepping scheme on velocity level even maintains physical consistency of the impulsive discretization. Hence, we propose a timestepping scheme based on Moreau's midpoint rule which enables to achieve not only compliance of the impact law but also of the non-penetration constraint. The choice of a decoupled and consecutive evaluation of the respective constraints can be interpreted as a not energy-consistent projection to the non-penetration constraint at the end of each time step. It is the implicit coupling of position and velocity level which yields satisfactory results. An implicit evaluation of the right hand side improves stability properties without additional cost. With the prox~function formulation, the overall set of nonsmooth equations is solved by a Newton scheme. Results from simulations of a slider-crank mechanism with unilateral constraints demonstrate the capability of our approach.

\section*{Keywords}
nonsmooth dynamics $\cdot$ timestepping scheme $\cdot$ Gear-Gupta-Leimkuhler method $\cdot$ unilateral contact $\cdot$ impact $\cdot$ slider-crank mechanism

\section{Introduction}\label{sec:introduction}
Dynamical motion with impacts plays an important role in the characterization of general mechanical systems at least after discretization in space. The monographs~\cite{Glo01,Aca08,Lei08,Pfe08,Ste11} summarize the state-of-the-art physical, mathematical and numerical setting of this kind of \emph{impacting mechanical systems}:
\begin{align}
  \vq\left(0\right)&=\vq_0\;,\label{eq:mechanical_system_initial_position}\\
  \vv\left(0\right)&=\vv_0\;,\label{eq:mechanical_system_initial_velocity}\\
  \dot{\vq} &= \vv\;,\label{eq:mechanical_system_position}\\
  \vM \dot{\vv} &= \vh + \vW^{T}\vlambda\;,\label{eq:mechanical_system_velocity}\\
  \vM \left(\vv_j^+-\vv_j^-\right) &= \vW_j^{T}\vLambda_j\;,\label{eq:mechanical_system_velocity_jump}\\
  0 \leq \vg\ &\bot \ \vlambda \geq 0\;,\label{eq:mechanical_system_contact}\\
  \text{if } \vg_j \leq 0 ,\text{ then } 0 \leq \dot{\vg}_j^+ +\vepsilon \dot{\vg}_j^-\ &\bot \ \vLambda_j \geq 0\;.\label{eq:mechanical_system_impact}
\end{align}
Starting from the initial conditions~\eqref{eq:mechanical_system_initial_position}-\eqref{eq:mechanical_system_initial_velocity}, the development of the system's state given by position~$\vq$ and velocity~$\vv$ is described by a non-impulsive behavior \eqref{eq:mechanical_system_velocity} almost everywhere. It is influenced by the generalized mass matrix~$\vM$ and right hand side forces~$\vh$. Due to the Signorini-Moreau condition~\eqref{eq:mechanical_system_contact}, closed or opening \emph{scleronomic} contact gaps~$\vg$ affect this type of motion by varying contact force parameters~$\vlambda$. The notation~$\vg\ \bot \ \vlambda$ stands for $\vg^T\vlambda=0$. The force parameters weight the columns of $\vW^T$ in the equations of motion \eqref{eq:mechanical_system_velocity}. For countable time instances $t_j$, the velocities jump enforced by an impact $\vLambda_j$ according to~\eqref{eq:mechanical_system_velocity_jump}. Newton's impact law~\eqref{eq:mechanical_system_impact} provides the respective relationship of active pre- and post-impact velocities using the kinematic coefficient of restitution~$\vepsilon$. With~\eqref{eq:mechanical_system_position}, the position $\vq$ can be calculated by the fundamental theorem of calculus for weakly differentiable functions.\par
\begin{figure}[ht]
  \centering
  \def\svgwidth{0.8\columnwidth}
  \import{}{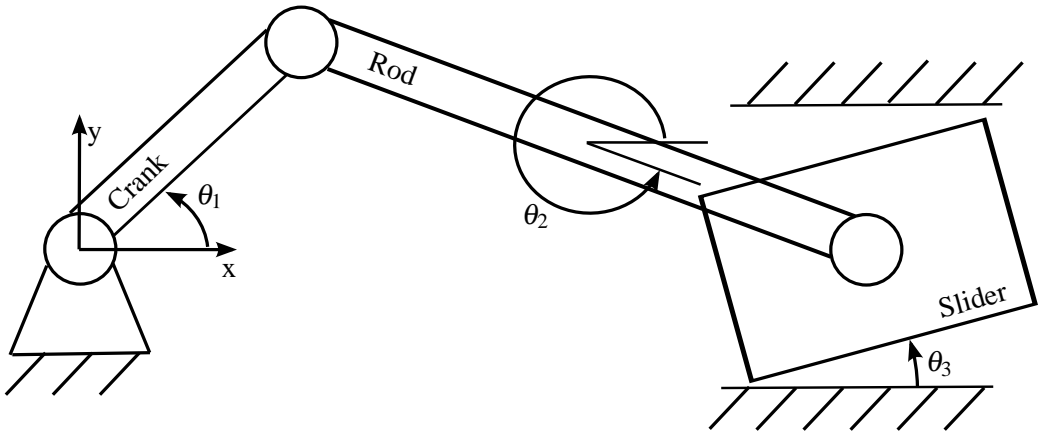}
  \caption{Slider-crank mechanism with unilateral constraints.}
  \label{fig:slider_crank_unilateral}
\end{figure}
One might think that the differential complementarity problem~\eqref{eq:mechanical_system_initial_position}-\eqref{eq:mechanical_system_impact} is integrated best by an event-driven time-integration strategy. However as event-driven schemes resolve the exact time of impact, they cannot consistently model Zeno phenomena, i.e. infinite impacts occurring in a finite time interval. Timestepping schemes discretize the equations of motion including the constraints consistently without resolving the exact transition points. Robustness benefits from this approach but the accuracy is comparably low.
\subsection{Moreau's midpoint rule}\label{sec_sub:moreau}
Because of the \emph{consistent approach} within \emph{timestepping method}s, we focus on a well-established representative, i.e. \emph{Moreau's midpoint rule}~\cite{Mor99}. It summarizes both impulsive and non-impulsive phases: first by calculating --in some sense-- the mean impulsive force within fixed time steps~$\Delta t$ and second by incorporating the results implicitly in a time-discretization on velocity level:
\begin{align}
  \vv_{n+1}&=\vv_n+\vM_M^{-1}\left(\vh_M\Delta t \ + \vW_M^T \vLambda_{n+1}\right)\;,\label{eq:dmoreau1}\\
  \vq_{n+1}&=\vq_n+\frac{\vv_{n+1}+\vv_{n}}{2} \Delta t \label{eq:dmoreau2}
\end{align}
with
\begin{align}
  \vM_M=\vM\left(\vq_n+\frac{\Delta t}{2}\vv_n\right)\;,\quad \vh_M=\vh\left(\vq_n+\frac{\Delta t}{2}\vv_n,\vv_n\right)\;,\quad \vW_M=\vW\left(\vq_n+\frac{\Delta t}{2}\vv_n\right)\;.\label{eq:dmoreau_evaluation}
\end{align}
Equation~\eqref{eq:dmoreau_evaluation} approximates the positions by their explicitly calculated midpoint values and the velocities by the respective explicit values at the beginning of the time step. The only unknown variables in \eqref{eq:dmoreau1}-\eqref{eq:dmoreau2} are the positions~$\vq_{n+1}$ and velocities~$\vv_{n+1}$ at the end of the time step and the discrete mean impulsive force $\vLambda_{n+1}$. With the explicit predictor 
\begin{align}
  \vg_{M}=\vg\left(\vq_n+\frac{\Delta t}{2}\vv_n\right)<0\;,\label{eq:activeset}
\end{align}
an active set of constraints indicated by the subscript $(\ )_{\text{red}}$ is determined and every corresponding interaction can be treated by Newton's impact law on velocity level. Next to classical impacts also closed contacts, opening contacts or impacts occurring without a normal contact impulse are naturally included. The complementarity formulation of Newton's impact law can be equivalently reformulated by dint of the prox~function related to a convex set $C\subset \MR$; this is often easier to solve numerically than a complementarity problem~\cite{Sch11a}. As the proximal point $\text{prox}_C(x)\in C$ of $x\in\MR$ is defined as 
\begin{align}
  \text{prox}_C(x)=\text{arg}\mathop{\text{min}}\limits_{x^* \in C} \left\|x-x^*\right\|\;,
\end{align} 
Newton's discrete impact law corresponds row-by-row to
\begin{align}
  \vLambda_{n+1,\text{red}}=\text{\textbf{prox}}_{\MR_0^+}\left(\vLambda_{n+1,\text{red}}-\vr\left(\dot{\vg}_{n+1,\text{red}}+\vepsilon\dot{\vg}_{n,\text{red}} \right) \right)\;.
\end{align} 
This expression can be solved iteratively with the --in the easiest case-- positive diagonal parameter matrix $\vr$ controlling the speed of convergence; in this work, we choose a Newton scheme without adapting $\vr$ as solution method. As termination criterion, the natural monotony test or a tolerance for the residuum can be employed.
\subsection{Problem statement}
Moreau's midpoint rule holds for general impacting mechanical systems. We reveal improvement possibilities at a glance of a planar impacting slider-crank mechanism~\cite{Flo10}.
\subsubsection{Slider-crank mechanism with unilateral constraints}
For the slider-crank mechanism in Figure~\ref{fig:slider_crank_unilateral}, angles $\vq=\left(\theta_1,\theta_2,\theta_3\right)^T$ and angular velocities $\vv=\left(\omega_1,\omega_2,\omega_3\right)^T$ rely on an absolute description concerning an inertial $x$-$y$-frame of reference. The crank (1) has mass~$m_1$, rotational inertia around the center of gravity~$J_1$ and length~$l_1$. The connecting rod (2) is similarly represented by $m_2$, $J_2$ and $l_2$. The slider (3) with $m_3$ and $J_3$ has height~$2b$ and length~$2a$. Its center of gravity $\left(x_3,y_3\right)$ is not fixed on one y-position but can move within a notch of height~$d$ and clearance~$c$. 
\begin{figure}[hbt]
  \centering
  \def\svgwidth{0.8\columnwidth}
  \import{}{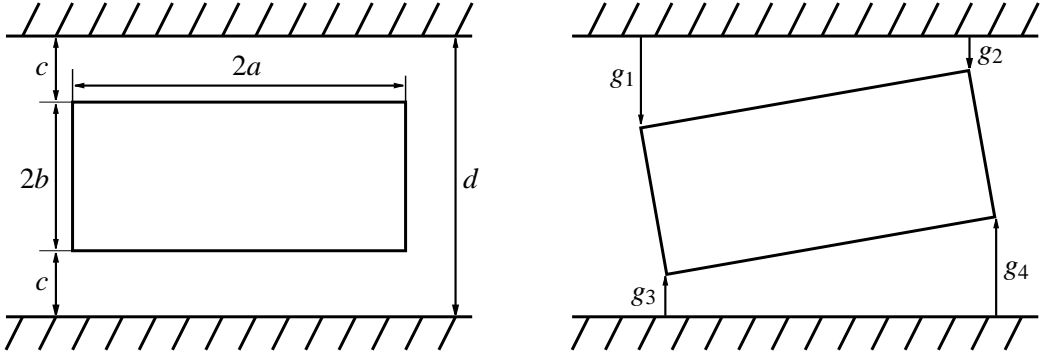}
  \caption{Definition of the gap functions for the slider-crank mechanism with unilateral constraints.}
  \label{fig:slider_crank_gaps}
\end{figure}
The gap functions are defined as illustrated in Figure~\ref{fig:slider_crank_gaps}:
\begin{align}
  g_1&=\frac{d}{2}-l_1\sin\theta_1-l_2\sin\theta_2+a\sin\theta_3-b\cos\theta_3\;,\\
  g_2&=\frac{d}{2}-l_1\sin\theta_1-l_2\sin\theta_2-a\sin\theta_3-b\cos\theta_3\;,\\
  g_3&=\frac{d}{2}+l_1\sin\theta_1+l_2\sin\theta_2-a\sin\theta_3-b\cos\theta_3\;,\\
  g_4&=\frac{d}{2}+l_1\sin\theta_1+l_2\sin\theta_2+a\sin\theta_3-b\cos\theta_3\;.
\end{align}
The tangential gap functions are neglected because the frictionless case is considered. Accordingly, the constraint matrix satisfies
\begin{align}
  \vW^T=\begin{pmatrix}
    -l_1\cos\theta_1 & -l_1\cos\theta_1 & l_1\cos\theta_1 & l_1\cos\theta_1 \\
    -l_2\cos\theta_2 & -l_2\cos\theta_2 & l_2\cos\theta_2 & l_2\cos\theta_2 \\
    a\cos\theta_3+b\sin\theta_3 & -a\cos\theta_3+b\sin\theta_3 & -a\cos\theta_3+b\sin\theta_3 & a\cos\theta_3+b\sin\theta_3 
  \end{pmatrix}\;.
\end{align}
Assuming gravitation $\vg$ in negative $y$-direction, the example fits exactly in the concept of a general impacting mechanical system~\eqref{eq:mechanical_system_initial_position}-\eqref{eq:mechanical_system_impact}:
\begin{align}
  \vM&=
  	\begin{pmatrix}
  	  J_1 + l_1^2 \left(\frac{m_1}{4}+m_2+m_3\right)& l_1 l_2 \cos \left(\theta_1-\theta_2\right)\left(\frac{m_2}{2}+m_3\right) & 0\\
  	  l_1 l_2 \cos \left(\theta_1-\theta_2\right)\left(\frac{m_2}{2}+m_3\right) & J_2 + l_2^2\left(\frac{m_2}{4}+m_3\right) & 0  \\
  	  0 & 0 & J_3
  	\end{pmatrix}\;,\\
  \vh&=
  	\begin{pmatrix}
  	  -l_1 l_2 \sin\left(\theta_1-\theta_2\right)\left(\frac{m_2}{2}+m_3\right) \omega_2^2 - g l_1 \cos \theta_1\left(\frac{m_1}{2}+m_2+m_3\right) \\
  	  l_1 l_2 \sin\left(\theta_1-\theta_2\right)\left(\frac{m_2}{2}+m_3\right) \omega_1^2 - g l_2 \cos \theta_2\left(\frac{m_2}{2}+m_3\right) \\ 
  	  0
  	\end{pmatrix}\;.
\end{align}
\subsubsection{Simulation results}
\begin{figure}[th]
  \centering
  \subfigure[$\vepsilon=0.1$]{\includegraphics[width=0.45\columnwidth]{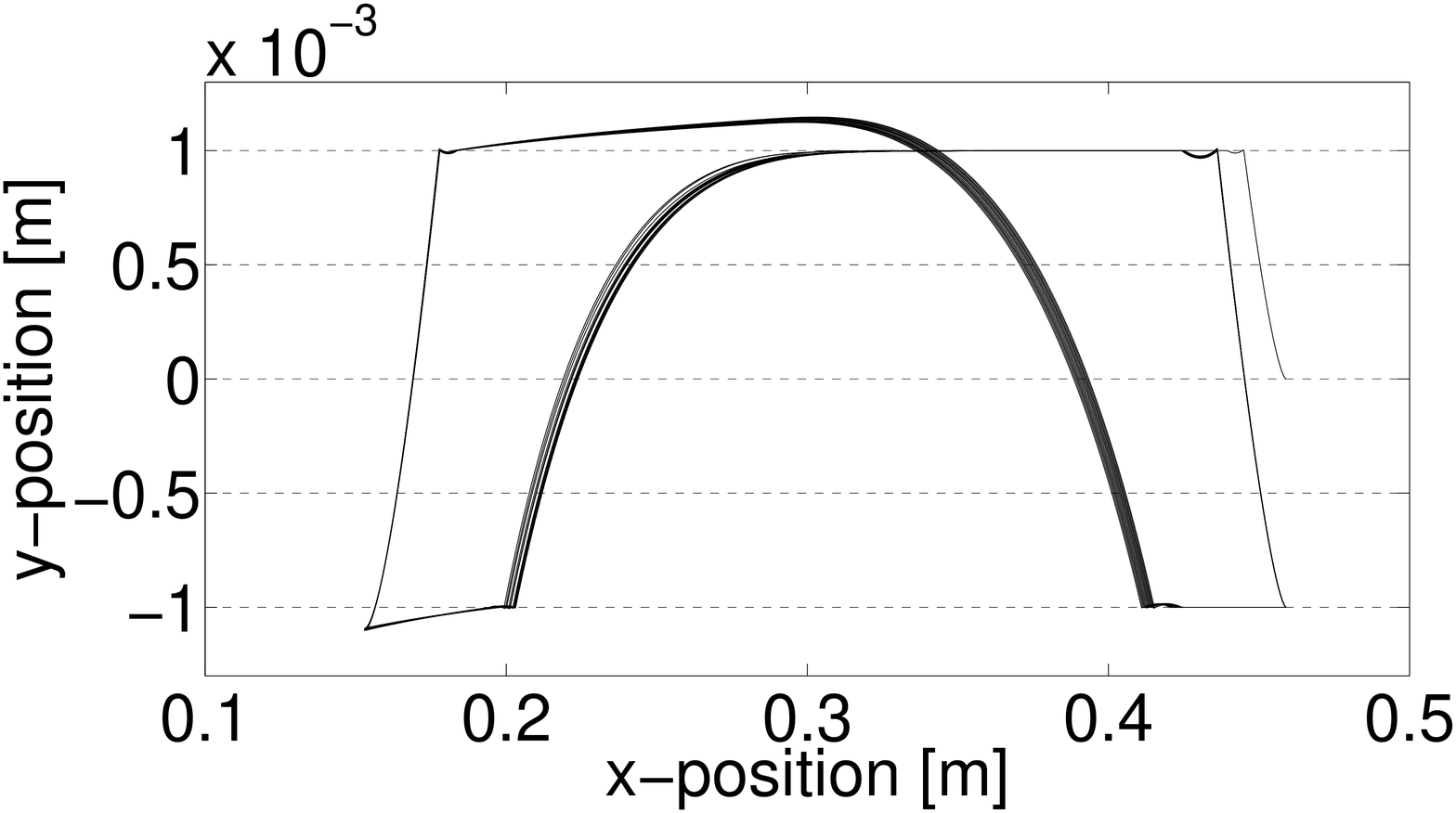}}
  \subfigure[$\vepsilon=0.4$]{\includegraphics[width=0.45\columnwidth]{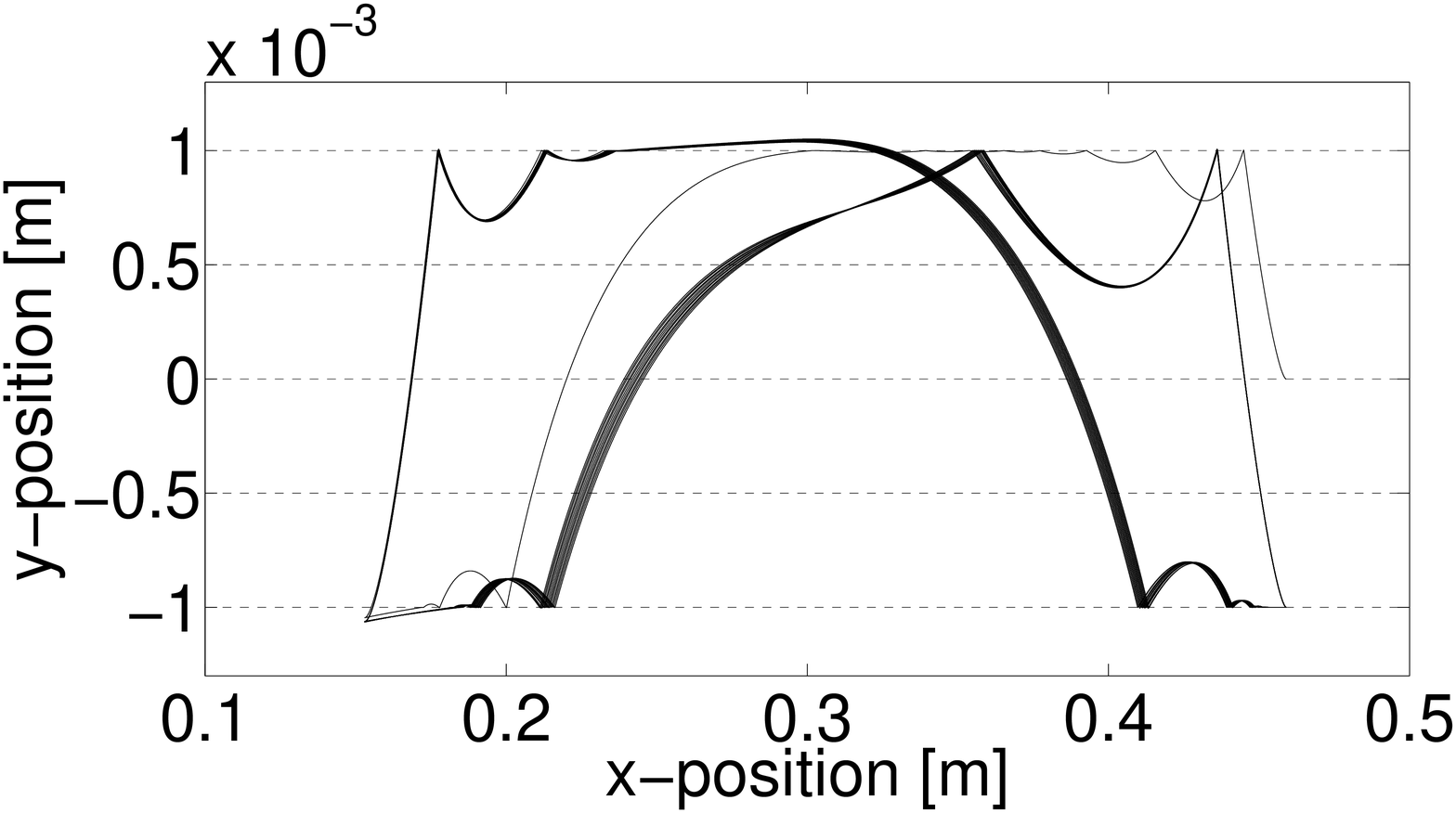}}\\
  \subfigure[$\vepsilon=0.6$]{\includegraphics[width=0.45\columnwidth]{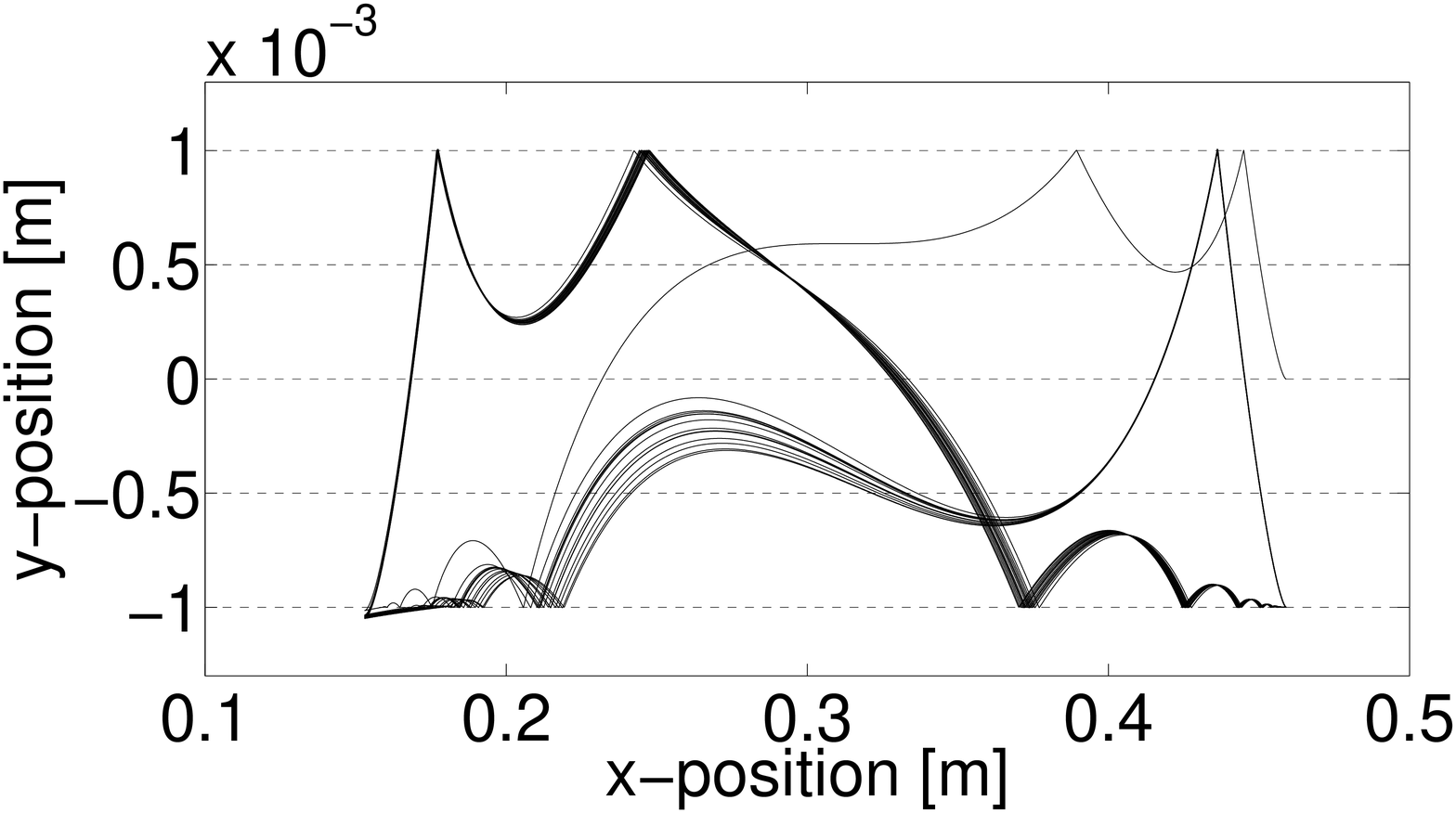}}
  \subfigure[$\vepsilon=0.9$]{\includegraphics[width=0.45\columnwidth]{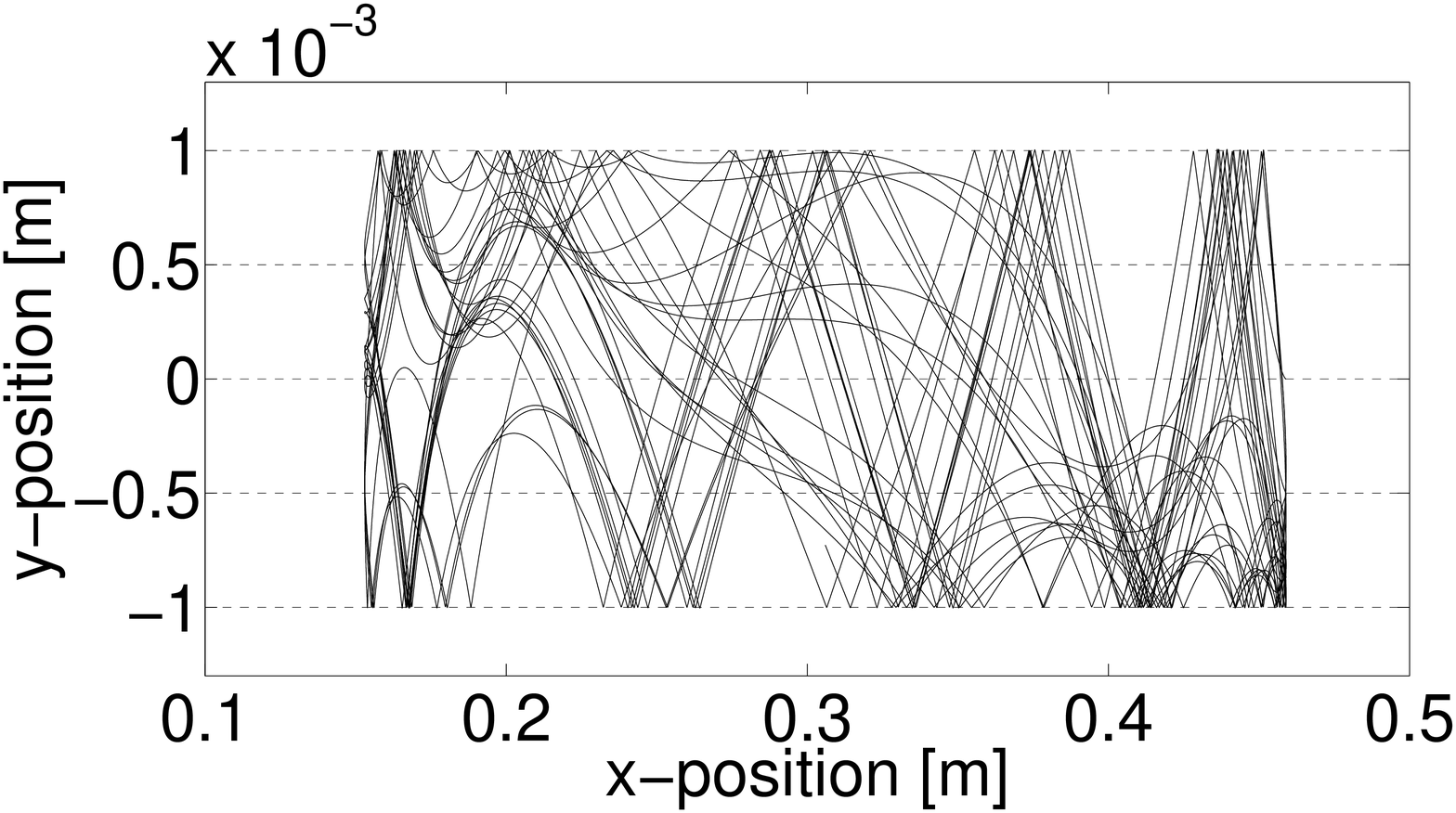}}
  \caption{Movement of the center of gravity of the slider (3) for different coefficients of restitution.}
  \label{fig:resmoreau}
\end{figure}

With the characteristics of~\cite{Flo10} reprinted in~Table~\ref{tab:slider_crank_characteristics}, we analyze the movement of the center of gravity of the slider (3) using Moreau's midpoint rule and the time step size $\Delta t=\unit[10^{-5}]{s}$. Figure~\ref{fig:resmoreau} shows the results for four different coefficients of restitution being the same for each contact possibility, respectively.
The curves presented here, as well as the graphs depicted in the original literature~\cite{Flo10}, show the violation of the non-penetration condition~\eqref{eq:mechanical_system_contact}, especially in the simulations with a low coefficient of restitution. In Figure~\ref{fig:resmoreau} (a), the center of gravity of the slider (3) exceeds $\unit[10^{-3}]{m}$ and in the later course it falls below the value $\unit[-10^{-3}]{m}$, what corresponds to a pervasion of the slider (3) with the bordering wall. Also for higher coefficients of restitution, the gap functions fall below zero, but of course for shorter time periods.\par
\begin{figure}[ht]
  \centering
  \includegraphics[width=1.0\columnwidth]{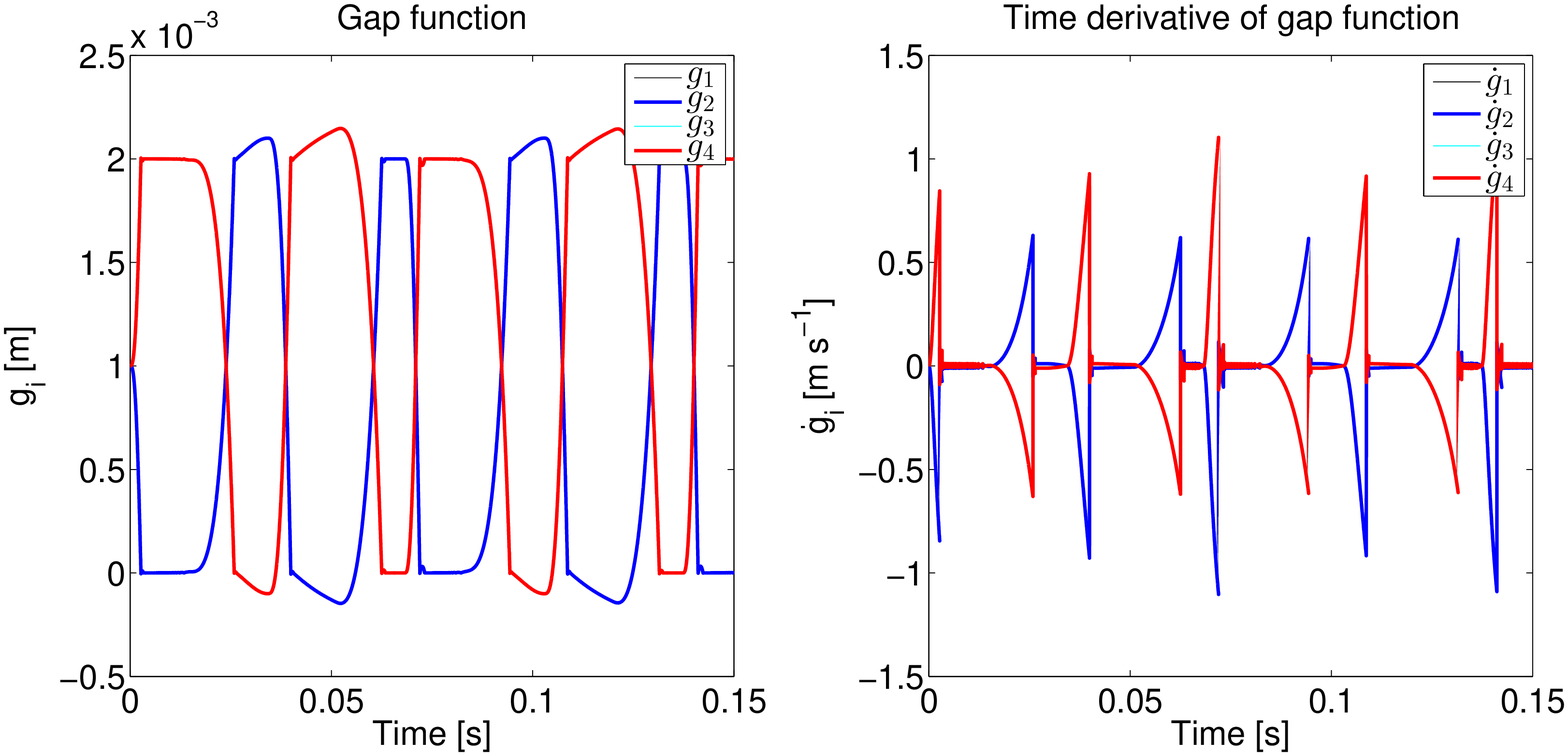}
  \caption{Gap functions and their time derivatives for $\epsilon=0.1$.}
  \label{fig:gapsmoreau}
\end{figure}

Figure \ref{fig:gapsmoreau} shows the time curve of the gap functions and their time derivatives for a coefficient of restitution $\vepsilon=0.1$.
The gap between the slider (3) and the surrounding wall is small as well as the initial configuration $\theta_{3_0}=0$ and the support at the center of gravity prevent the revolution of the slider (3). Hence, pairs of gap functions on opposite sides appear almost symmetric. Obviously, the non-penetration condition is violated. A detailed view shows the \emph{drift-off~effect}: the gap drifts approximately linearly to negative values while the gap velocity is slightly negative. However, drift does not have a dominant effect for this configuration because the negative gap functions remain comparatively small in contrast to the geometric dimensions. When contacts stay closed for longer time periods, the drift-off effect will not be negligible anymore.
\begin{table}[hbt]
  \centering
  \begin{tabular}{ll}
    \hline
    Geometrical characteristics & $l_1=\unit[0.1530]{m}$\\
    \ & $l_2=\unit[0.3060]{m}$\\
    \ & $a=\unit[0.0500]{m}$\\
    \ & $b=\unit[0.0250]{m}$\\
    \ & $c=\unit[0.0010]{m}$\\ \hline
    Inertia properties & $m_1=\unit[0.0380]{kg}$\\
    \ & $m_2=\unit[0.0380]{kg}$\\
    \ & $m_3=\unit[0.0760]{kg}$\\
    \ & $J_1=\unit[7.4\cdot10^{-5}]{kgm^2}$\\
    \ & $J_2=\unit[5.9\cdot10^{-4}]{kgm^2}$\\
    \ & $J_3=\unit[2.7\cdot10^{-6}]{kgm^2}$\\ \hline
    \ Force elements & $g=\unit[9.81]{m/s^2}$\\ \hline
    \ Initial conditions & $\theta_{10}=0.0$\\
    \ & $\theta_{2_0}=0.0$\\
    \ & $\theta_{3_0}=0.0$\\
    \ & $\omega_{1_0}=\unit[150.0]{1/s}$\\
    \ & $\omega_{2_0}=\unit[-75.0]{1/s}$\\
    \ & $\omega_{3_0}=\unit[0.0]{1/s}$\\ \hline
  \end{tabular}
  \caption{Characteristics of the slider-crank mechanism with unilateral constraints.}
  \label{tab:slider_crank_characteristics}
\end{table}
\subsubsection{Outline}
For systems with only persistent contacts, one could consider the constraints on position level. For systems which undergo impacts in addition, this would yield non-consistent discretizations. Hence, when both impacts \emph{and} longer periods of closed contacts occur, neither of these two approaches, i.e. \emph{neither} on position \emph{nor} velocity level seems to be satisfactory. Gear, Gupta and Leimkuhler~\cite{Gea85} proposed a solution to a related problem for only persistent contacts: they enforce constraints on position and velocity level at the same time. The additional constraint equation is compensated by a second set of Lagrange multipliers. The purpose of this work is to apply the \emph{Gear-Gupta-Leimkuhler method} to systems with \emph{unilateral} constraints to overcome the drift-off~effect for closed contacts as roughly indicated in~\cite{Aca11b}. Thereby, we summarize and extent our student work~\cite{Sch12c}. First, we present the Gear-Gupta-Leimkuhler method for a slider-crank mechanism without clearance. The application to unilateral contacts shows that a decoupled strategy satisfying velocity and position level constraints one after the other is not energy-consistent. A unified formulation which takes into account the impact law as well as the non-penetration constraint at the same time turns out to be a successful approach. We close the paper with some open questions for future research directions.

\section{Gear-Gupta-Leimkuhler method for persistent contacts}
We explain the expected effect of the Gear-Gupta-Leimkuhler method~\cite{Gea85} on the numerical solution of a unilaterally constrained mechanical system with the help of a bilaterally constrained slider-crank mechanism adapted from~\cite{Flo10}. The drift-off effect is analyzed for constraint formulations on position, velocity and acceleration level as well as for the Gear-Gupta-Leimkuhler formulation~\cite{Fueh91}.
\subsection{Slider-crank mechanism with bilateral constraints}
\begin{figure}[hbt]
  \centering
  \def\svgwidth{0.8\columnwidth}
  \import{}{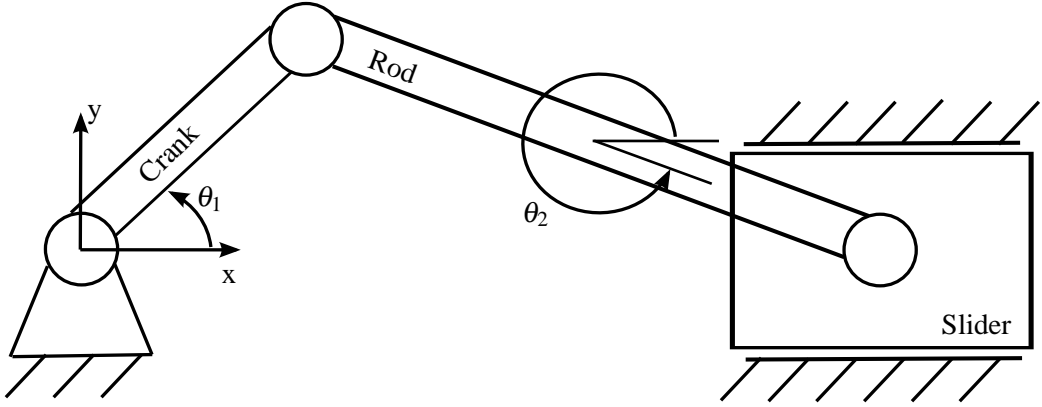}
  \caption{Slider-crank mechanism with bilateral constraints.}
  \label{fig:slider_crank_bilateral}
\end{figure}

Figure~\ref{fig:slider_crank_bilateral} shows a bilaterally constrained slider-crank mechanism.
The angles $\vq=\left(\theta_1,\theta_2\right)^T$ are chosen as generalized coordinates, the angular velocities $\vv=\left(\omega_1,\omega_2\right)^T$ as generalized velocities. With the same notation as in Section~\ref{sec:introduction}, we gain equations of motion which are more specific than stated in~\eqref{eq:mechanical_system_initial_position}-\eqref{eq:mechanical_system_impact}: unilateral contacts condense to bilateral constraints and impacts never occur. The generalized mass matrix satisfies
\begin{align}
\vM=
	\begin{pmatrix}
	  J_1 + l_1^2 \left(\frac{m_1}{4}+m_2+m_3\right)& l_1 l_2 \cos \left(\theta_1-\theta_2\right)\left(\frac{m_2}{2}+m_3\right)\\
	  l_1 l_2 \cos \left(\theta_1-\theta_2\right)\left(\frac{m_2}{2}+m_3\right) & J_2 + l_2^2\left(\frac{m_2}{4}+m_3\right)
	\end{pmatrix}
\end{align}
and the vector of generalized forces is given by
\begin{align}
\vh=
	\begin{pmatrix}
	  -l_1 l_2 \sin\left(\theta_1-\theta_2\right)\left(\frac{m_2}{2}+m_3\right) \omega_2^2 - g l_1 \cos \theta_1\left(\frac{m_1}{2}+m_2+m_3\right) \\
	  l_1 l_2 \sin\left(\theta_1-\theta_2\right)\left(\frac{m_2}{2}+m_3\right) \omega_1^2 - g l_2 \cos \theta_2\left(\frac{m_2}{2}+m_3\right)
	\end{pmatrix}\;.
\end{align}
The bilateral constraint holds the slider (3) at a fixed y-position
\begin{align}
  g=l_1 \sin \theta_1 + l_2 \sin \theta_2=0
\end{align}
by causing a constraint force in direction of 
\begin{align}
  \vW^T=\begin{pmatrix}
    l_1 \cos \theta_1\\
    l_2 \cos \theta_2
  \end{pmatrix}\;.
\end{align}
\subsection{Simulation results}
\begin{figure}[hbt]
  \centering
  \includegraphics[width=0.8\columnwidth]{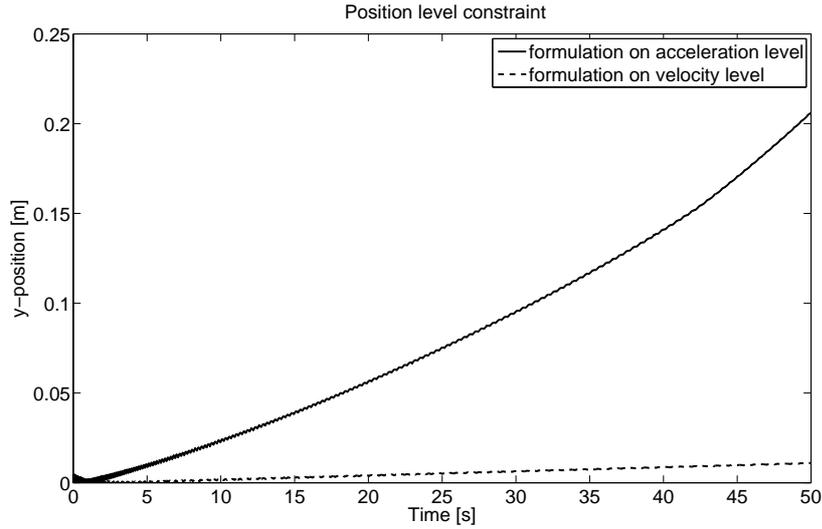}
  \caption{Drift-off effect of the slider (3) for equations formulated on velocity and acceleration level.}
  \label{fig:slider_crank_bilateral_drift}
\end{figure}

The simulations are accomplished with the characteristics of Table~\ref{tab:slider_crank_characteristics} and the time step size $\Delta t=\unit[10^{-4}]{s}$. A direct computation of the constraint compliance considering the constraint $g$ on position level yields a differential algebraic system of index 3. It is known to be badly conditioned and e.g. scaling of the constraint equation yields an heuristic improvement concerning the stability of the numerical integration scheme~\cite{Hai10}. Instead of that, we focus on replacing the constraint by its respective time derivatives, which improves the robustness of numerical solvers consistently from an analytic point of view. Arnold~\cite{Arn09} mentions this strategy in the context of \emph{index reduction}. A draw back of index reduction is the drift-off effect~\cite{Arn09}. Figure~\ref{fig:slider_crank_bilateral_drift} shows the roughly linear development of the y-position of the slider (3) for a long-time simulation of the index 2 system, i.e. considering the constraint $\dot{g}$ on velocity level. 
Figure~\ref{fig:slider_crank_bilateral_drift} also displays the drift-off effect for the constraint formulation on acceleration level~$\ddot{g}$, i.e. for the index 1 system. As presented in~\cite{Fueh91}, the drift-off is expected to be parabolic and in fact the y-position increases with rising gradient. As a compromise of both robust simulation and asymptotically little drift-off, one usually considers the constraints on velocity level. To even overcome the linear drift-off effect in the index 2 system, Gear, Gupta and Leimkuhler proposed a formulation which considers the constraints on position \emph{and} velocity level simultaneously~\cite{Gea85}. The original index 2 system extends to 
\begin{align}
  	  \dot{\vq} &= \vv\color{blue}{+\vW^{T}\vpsi}\;,\\
      \vM\vv &= \vh+\vW^{T}\vlambda\;,\\
  	  \dot{\vg} &= \vnull\;,\\ 
  	  \color{blue}{\vg} &= \color{blue}{\vnull}\;.
\end{align}
The Lagrange multiplier $\vpsi$ compensates the added equation and the constraint is satisfied on position as well as on velocity level maintaining the stability of the index 2 formulation. 

\section{Gear-Gupta-Leimkuhler method for unilateral contacts}
We analyze two extensions of Moreau's midpoint rule (cf.~Section~\ref{sec_sub:moreau}). A decoupled approach turns out not to be energy-consistent. A unified approach meets our expectations but demands the computational effort of an implicit solution scheme. 
\subsection{Decoupled approach}
The adaption of Moreau's midpoint rule is performed by adding a correction term enforcing the non-penetration constraint at the end of each time step:  
\begin{align}
  \vv_{n+1}&=\vv_n+\vM_M^{-1}\left(\vh_M\Delta t \ + \vW_M^T \vLambda_{n+1}\right)\;,\label{eq:ggl1}\\
  \vLambda_{n+1,\text{red}}&=\text{\textbf{prox}}_{\MR_0^+}\left(\vLambda_{n+1,\text{red}}-\vr\left(\dot{\vg}_{n+1,\text{red}}+\vepsilon\dot{\vg}_{n,\text{red}} \right) \right)\;,\label{eq:prox_velocity}\\
  \vq_{n+1}&=\vq_n+\frac{\vv_{n+1}+\vv_{n}}{2} \Delta t\ + \vW_M^T \vPsi_{n+1}\;,\label{eq:ggl2}\\
  \vPsi_{n+1}&=\text{\textbf{prox}}_{\MR_0^+}\left(\vPsi_{n+1}-\vr\vg_{n+1}\right)\;.\label{eq:prox_position}
\end{align}
As in Moreau's midpoint rule, the calculation of the velocities~$\vv_{n+1}$ is achieved by using the average Lagrange multiplier $\vLambda_{n+1}$. This computation is decoupled from the calculation of the positions $\vq_{n+1}$ with the average Lagrange multiplier $\vPsi_{n+1}$. Consecutively, the velocities~$\vv_{n+1}$ are used to determine an explicit forecast:
\begin{align}
  \bar{\vq}_{n+1}=\vq_n+\frac{\vv_{n+1}+\vv_n}{2}\Delta t\;.
\end{align}
It is used as an initial value for the iterative computation of $\vq_{n+1}$ and $\vPsi_{n+1}$.\par
The positive effect of low computational effort, due to a decoupled calculation of the two different vectors of Lagrange multipliers, is subtended by the low physical accuracy of the results. This is clarified by the development of the entire energy content of the slider-crank mechanism with unilateral constraints using the characteristics of Table~\ref{tab:slider_crank_characteristics} (cf.~Figure~\ref{fig:energycontent}). No energy sources are applied but for a coefficient of restitution $\epsilon=0.1$, our simulation results with time step size $\Delta t=\unit[10^{-5}]{s}$ reveal a fluctuating entire energy content and do not show the expected decreasing trend. Hence, the timestepping scheme~\eqref{eq:ggl1}-\eqref{eq:prox_position} does not provide a valid and physical accurate model of a system underlying unilateral constraints.
\subsection{Unified approach}
What is the problem in~\eqref{eq:ggl1}-\eqref{eq:prox_position}? The equations of motion are derived via an energy principle, the impact law results from Newton's admittedly kinematic considerations. However, the additional term $\vW^T_M\vPsi_{n+1}$ does not correspond to any physical principle but can be interpreted as part of the Karush-Kuhn-Tucker conditions for a \emph{projection at the end of each time step} concerning the Euclidean metric:
\begin{align}
  \min_{\vq_{n+1}}\norm{\vq_{n+1}-\bar{\vq}_{n+1}}^2\;,\\
  \vg\left(\vq_{n+1}\right)\geq\vnull\;.
\end{align}
It is not astonishing that the entire energy content oscillates. Studer~\cite{Stu09} mentions this type of discretization discussing the bouncing ball example. As a workaround, it is suggested to introduce a penetration tolerance depending on the specific setting. To our opinion, the term $\vW^T_M\vPsi_{n+1}$ has to be coupled with~\eqref{eq:ggl1} to ensure a physical accurate behavior in general. Our proposition is presented in the next section.
\subsubsection{Discretization scheme}
\begin{figure}
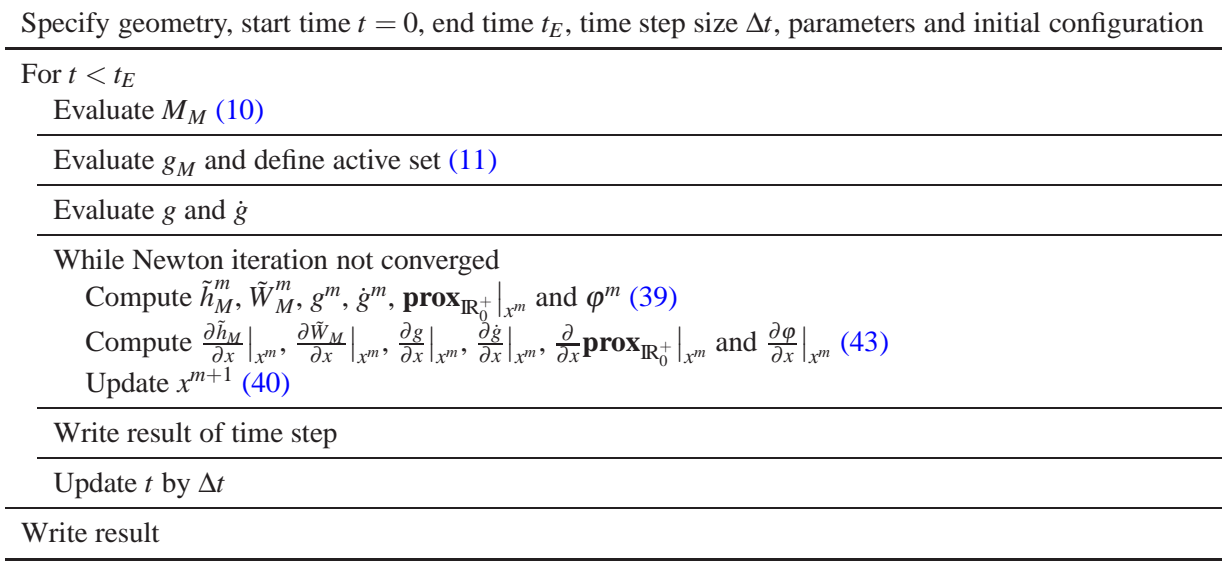

  \centering
  \begin{tabular}{lll}
  \toprule
  \multicolumn{3}{l}{Specify geometry, start time $t=0$, end time $t_E$, time step size $\Delta t$, parameters and initial configuration}\\
  \midrule 
  \multicolumn{3}{l}{For $t<t_E$}\\
   & \multicolumn{2}{l}{Evaluate $\vM_M$ \textcolor{blue}{\eqref{eq:dmoreau_evaluation}}}\\
  \cmidrule{2-3}
   & \multicolumn{2}{l}{Evaluate $\vg_M$ and define active set \textcolor{blue}{(\ref{eq:activeset})}}\\
  \cmidrule{2-3}
   & \multicolumn{2}{l}{Evaluate $\vg$ and $\dot{\vg}$}\\
  \cmidrule{2-3}
   & \multicolumn{2}{l}{While Newton iteration not converged}\\
   & & Compute $\tilde{\vh}_{M}^m$, $\tilde{\vW}_{M}^m$, $\vg^m$, $\dot{\vg}^m$, $\text{\textbf{prox}}_{\MR^+_0}\big|_{\vx^m}$ and $\vvarphi^m$ \textcolor{blue}{\eqref{eq:magnime2}}\\
   & & Compute $\frac{\partial \tilde{\vh}_M}{\partial \vx}\big|_{\vx^m}$, $\frac{\partial \tilde{\vW}_M}{\partial \vx}\big|_{\vx^m}$, $\frac{\partial \vg}{\partial \vx}\big|_{\vx^m}$, $\frac{\partial \dot{\vg}}{\partial \vx}\big|_{\vx^m}$, $\frac{\partial}{\partial \vx}\text{\textbf{prox}}_{\MR^+_0}\big|_{\vx^m}$ and $\frac{\partial \vvarphi}{\partial\vx}\big|_{\vx^m}$ \textcolor{blue}{\eqref{eq:magnime3}}\\
   & & Update $\vx^{m+1}$ \textcolor{blue}{\eqref{eq:newton_iteration}}\\
  \cmidrule{2-3}
   & \multicolumn{2}{l}{Write result of time step}\\
  \cmidrule{2-3}
  & \multicolumn{2}{l}{Update $t$ by $\Delta t$}\\
  \midrule
  \multicolumn{3}{l}{Write result}\\
  \bottomrule
  \end{tabular}\caption{Flowchart of the proposed unified timestepping scheme.}
  \label{fig:flow}
\end{figure}

A possible coupling is the implicit evaluation of the constraint matrix~$\tilde{\vW}_M=\vW\left(\frac{\vq_{n+1}+\vq_n}{2}\right)$ maintaining the nice properties of the midpoint concept, e.g. symplecticity~\cite{Hai06}. As this strategy already enforces the solution of a nonlinear system of equations, we additionally evaluate the generalized force vector $\tilde{\vh}_M=\vh\left(\frac{\vq_{n+1}+\vq_n}{2},\frac{\vv_{n+1}+\vv_n}{2}\right)$ implicitly to benefit from a more stable discretization of important stiffness contributions. The generalized mass matrix comprises geometric nonlinearities and its implicit evaluation needs comparatively large effort. Hence, an explicit evaluation $\vM_M$ is chosen:
\begin{align}
  \vq_{n+1}&=\vq_n+\frac{\vv_{n+1}+\vv_{n}}{2} \Delta t\ + \vW^T\left(\frac{\vq_{n+1}+\vq_n}{2}\right) \vPsi_{n+1}\label{eq:gglfr2}\;,\\
  \vv_{n+1}&=\vv_n+\vM_M^{-1}\left[\vh\left(\frac{\vq_{n+1}+\vq_n}{2},\frac{\vv_{n+1}+\vv_n}{2}\right)\Delta t \ + \vW^T\left(\frac{\vq_{n+1}+\vq_n}{2}\right) \vLambda_{n+1}\right]\;.\label{eq:gglfr1}
\end{align}
Adding the active constraints and defining a nonlinear system of equations, only the dependency on the unknown variables 
\begin{align}
  \vx_{n+1,\mathrm{red}}=\left(\vq_{n+1}^T\ \vv_{n+1}^T\ \vLambda^T_{n+1,\mathrm{red}}\ \vPsi^T_{n+1,\mathrm{red}}\right)^T
\end{align}
is interesting:
\begin{align}
  \vvarphi_{\mathrm{red}}\left(\vx_{n+1,\mathrm{red}}\right)=\begin{pmatrix}
  \vq_{n+1}-\vq_n-\frac{\vv_{n+1}+\vv_{n}}{2} \Delta t\ - \tilde{\vW}^T_M\left(\vq_{n+1}\right) \vPsi_{n+1,\mathrm{red}}\\
  \vv_{n+1}-\vv_n-\vM_M^{-1}\left[\tilde{\vh}_M\left(\vq_{n+1},\vv_{n+1}\right)\Delta t \ + \tilde{\vW}^T_M\left(\vq_{n+1}\right) \vLambda_{n+1,\mathrm{red}}\right]\\
  \vLambda_{n+1,\mathrm{red}}-\text{\textbf{prox}}_{\MR^+_0}\left(\vLambda_{n+1,\mathrm{red}}-\vr\left(\dot{\vg}_{n+1,\mathrm{red}}+\vepsilon \dot{\vg}_{n,\mathrm{red}}\right)\right)\\
  \vPsi_{n+1,\mathrm{red}}-\text{\textbf{prox}}_{\MR^+_0}\left(\vPsi_{n+1,\mathrm{red}}-\vr \vg_{n+1,\mathrm{red}} \right)
  \end{pmatrix}=\vnull\;.\label{eq:magnime2}
\end{align}
In contrast to Section~\ref{sec:introduction}, the discretizations $\tilde{\vh}_M=\tilde{\vh}_M\left(\vq_{n+1},\vv_{n+1}\right)$ and $\tilde{\vW}_M=\tilde{\vW}_M\left(\vq_{n+1}\right)$ explicitly depend on the unknown values $\vq_{n+1}$ and $\vv_{n+1}$. The roots of the reduced system of equations $\vvarphi_{\mathrm{red}}\left(\vx_{\mathrm{red}}\right)$ can be solved using Newton's method
\begin{align}
  \vx_{n,\mathrm{red}}^{m+1}=\vx_{n,\mathrm{red}}^m- \left(\frac{\partial \vvarphi_{\mathrm{red}}}{\partial \vx_{\mathrm{red}}} \biggr|_{\vx_{n,\mathrm{red}}^m}\right)^{-1} \vvarphi_{\mathrm{red}}\left(\vx_{n,\mathrm{red}}^m\right)\;.\label{eq:newton_iteration}
\end{align}
With the scleronomic gap functions being simplified concerning effective evaluations
\begin{align}
  \vg_{n+1}&=W_{n+1}q_{n+1}\approx\vg_n+\tilde{\vW}_M\left(\vq_{n+1}\right)\frac{\vv_{n+1}+\vv_{n}}{2} \Delta t + \tilde{\vW}_M\left(\vq_{n+1}\right)\tilde{\vW}_M^T\left(\vq_{n+1}\right) \vPsi_{n+1}\label{eq:gglfr2g}\;,\\
  \dot{\vg}_{n+1}&=W_{n+1}v_{n+1}\approx\dot{\vg}_n+\tilde{\vW}_M\left(\vq_{n+1}\right)\vM_M^{-1}\left[\tilde{\vh}_M\left(\vq_{n+1},\vv_{n+1}\right)\Delta t + \tilde{\vW}^T_M\left(\vq_{n+1}\right) \vLambda_{n+1}\right]\;,\label{eq:gglfr1g}
\end{align}
the derivative of $\vvarphi_{\mathrm{red}}\left(\vx_{\mathrm{red}}\right)$ with respect to $\vx_{\mathrm{red}}$ can be deduced by eliminating rows and columns corresponding to inactive contacts from the following matrix
\begin{align}
  \frac{\partial \vvarphi}{\partial \vx}\biggr|_{\vx_n^m}=\left.\left(\begin{tabular}{cccc}
  $\vI-\frac{\partial\tilde{\vW}_M^T}{\partial\vq}\vPsi$ & $-\vI\frac{\Delta t}{2}$ & $\vnull$ & $-\tilde{\vW}_M^T$ \\
  $-\vM_M^{-1}\left(\frac{\partial\tilde{\vh}_M}{\partial \vq}\Delta t+\frac{\partial\tilde{\vW}_M^T}{\partial\vq}\vLambda \right)$ & $\vI-\vM_M^{-1}\frac{\partial\tilde{\vh}_M}{\partial \vv}\Delta t$ & $-\vM_M^{-1}\tilde{\vW}_M^T$ & $\vnull$\\
 \multicolumn{4}{c}{$\frac{\partial}{\partial \vx} \left( \vLambda_{n+1,\mathrm{red}}-\text{\textbf{prox}}_{\MR^+_0}\left(\vLambda_{n+1,\mathrm{red}}-\vr\left(\dot{\vg}_{n+1,\mathrm{red}}+\vepsilon \dot{\vg}_{n,\mathrm{red}}\right)\right)\right)$}\\
 \multicolumn{4}{c}{$\frac{\partial}{\partial \vx}\left( \vPsi_{n+1,\mathrm{red}}-\text{\textbf{prox}}_{\MR^+_0}\left(\vPsi_{n+1,\mathrm{red}}-\vr \vg_{n+1,\mathrm{red}} \right)\right)$}
  \end{tabular}\right)\right|_{\vx_n^m}\;.\label{eq:magnime3}
\end{align}
With an appropriate function $\vf$, a distinction of cases is necessary for the rows containing the prox~function:
\begin{align}
\frac{\partial}{\partial \vx}\left(\text{\textbf{prox}}_{\MR^+_0}\left(\vf(\vx)\right)\right)\biggr|_{\vx_n^m}=\begin{cases}
   \frac{\partial}{\partial \vx}f(\vx) \big|_{\vx_n^m},  & \text{if }f(\vx_n^m)>\vnull\\
   0, & \text{else }
 \end{cases}\;.
\end{align}
Analyzing exemplary the more difficult \texttt{if}-case of the prox~function, the derivatives of the third row of $\vvarphi$ are
\begin{align}
  \frac{\partial \vvarphi_3}{\partial \vq}&=\vr\left( \tilde{\vW}_M \vM^{-1}_M \left(\frac{\partial \tilde{\vh}_M}{\partial \vq}\Delta t + \frac{\partial \tilde{\vW}_M^T}{\partial \vq} \vLambda \right) + \frac{\partial \tilde{\vW}_M}{\partial \vq}\vM^{-1}_M \left( \tilde{\vh}_M \Delta t + \tilde{\vW}_M^T \vLambda \right)\right)\;,\\
  \frac{\partial \vvarphi_3}{\partial \vv}&=\vr\tilde{\vW}_M \vM^{-1}_M \frac{\partial \tilde{\vh}_M}{\partial \vv}\Delta t\;,\\
  \frac{\partial \vvarphi_3}{\partial \vLambda}&=\vr\tilde{\vW}_M \vM^{-1}_M \tilde{\vW}_M^T\;,\\
  \frac{\partial \vvarphi_3}{\partial \vPsi}&=\vnull\;.
\end{align}
For the fourth row, it is 
\begin{align}
  \frac{\partial \vvarphi_4}{\partial \vq}&=\vr\left(\frac{\partial \tilde{\vW}_M}{\partial \vq} \frac{\vv_n^m+\vv_n}{2} \Delta t + \frac{\partial \tilde{\vW}_M}{\partial \vq} \tilde{\vW}_M^T \vPsi+ \tilde{\vW}_M \frac{\partial \tilde{\vW}_M^T}{\partial \vq} \vPsi\right)\;,\\
  \frac{\partial \vvarphi_4}{\partial \vv}&=\vr\tilde{\vW}_M \frac{\Delta t }{2}\;,\\
  \frac{\partial \vvarphi_4}{\partial \vLambda}&=\vnull\;,\\
  \frac{\partial \vvarphi_4}{\partial \vPsi}&=\vr\tilde{\vW}_M\tilde{\vW}_M^T\;.
\end{align}
As long as the active set is empty, all Lagrange multipliers are equal to zero and the system of equations is solved without the constraint part. As soon as the active set is not empty, the system of equations is extended by Newton's impact law and the non-penetration constraint, respectively. The algorithm is set up as shown in Figure~\ref{fig:flow}.
\subsubsection{Simulation results}
\begin{figure}[hbt]
  \centering
  \subfigure[$\epsilon=0.1$]{\includegraphics[width=0.45\columnwidth]{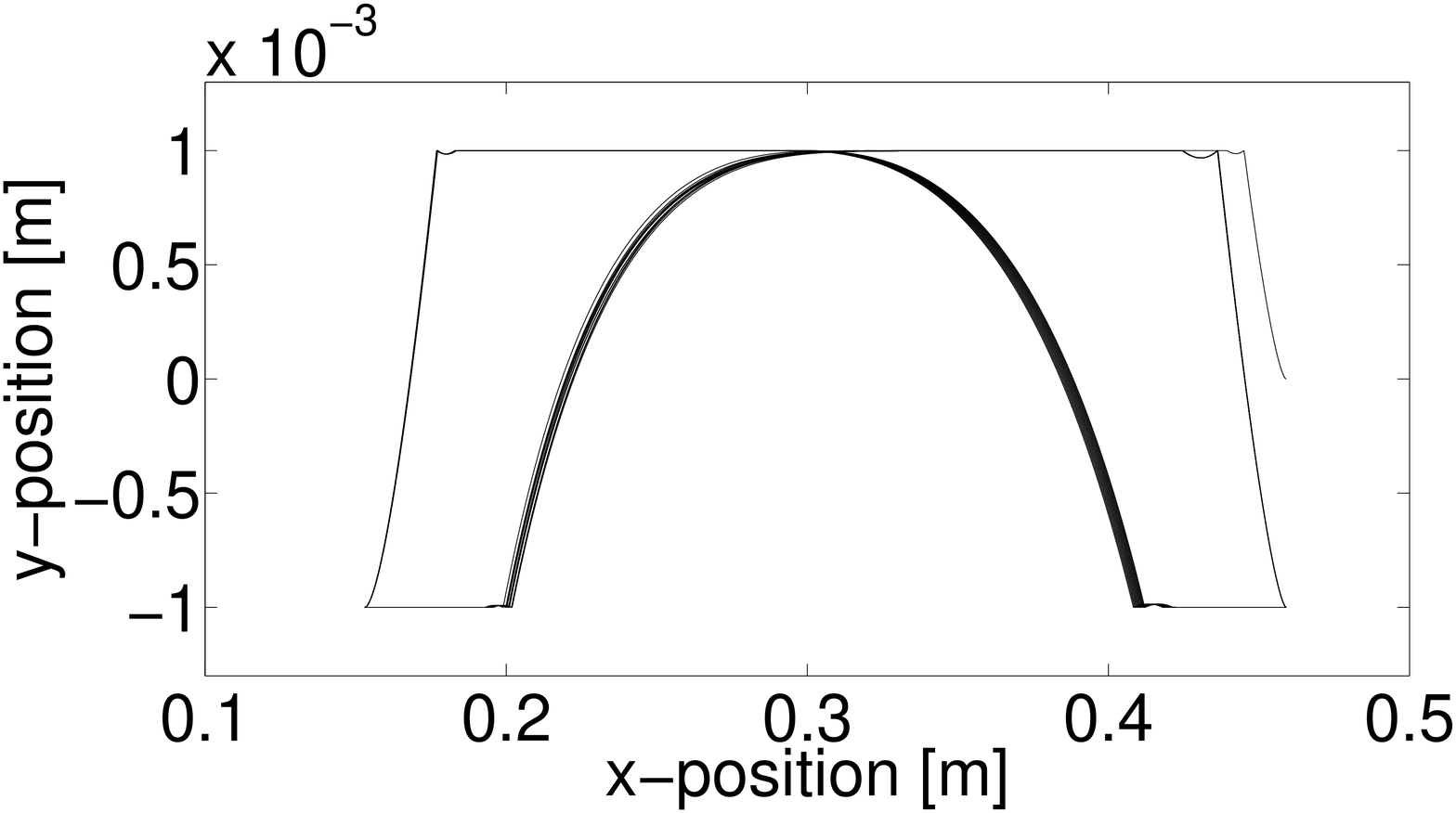}}
  \subfigure[$\epsilon=0.4$]{\includegraphics[width=0.45\columnwidth]{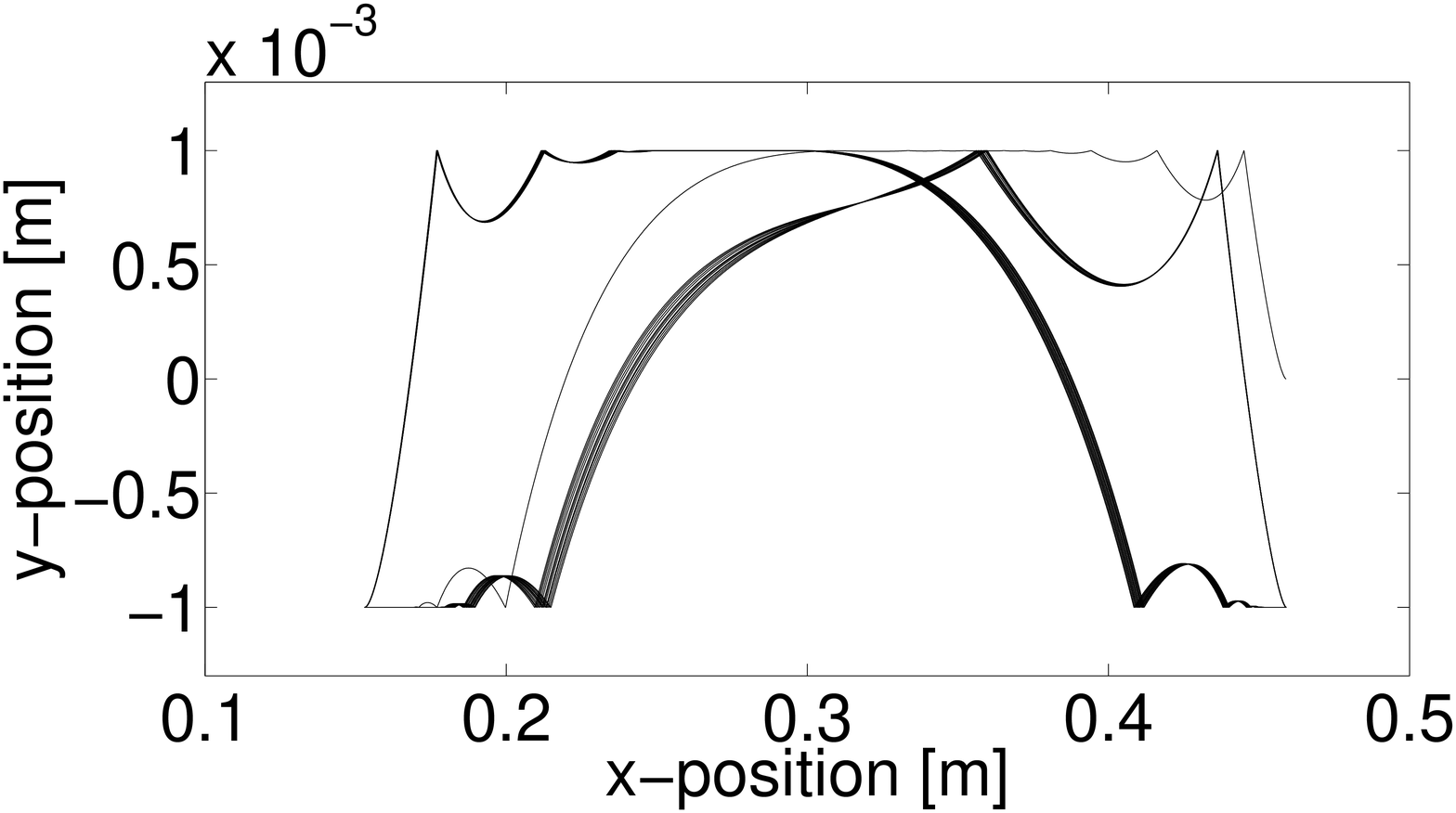}}\\
  \subfigure[$\epsilon=0.6$]{\includegraphics[width=0.45\columnwidth]{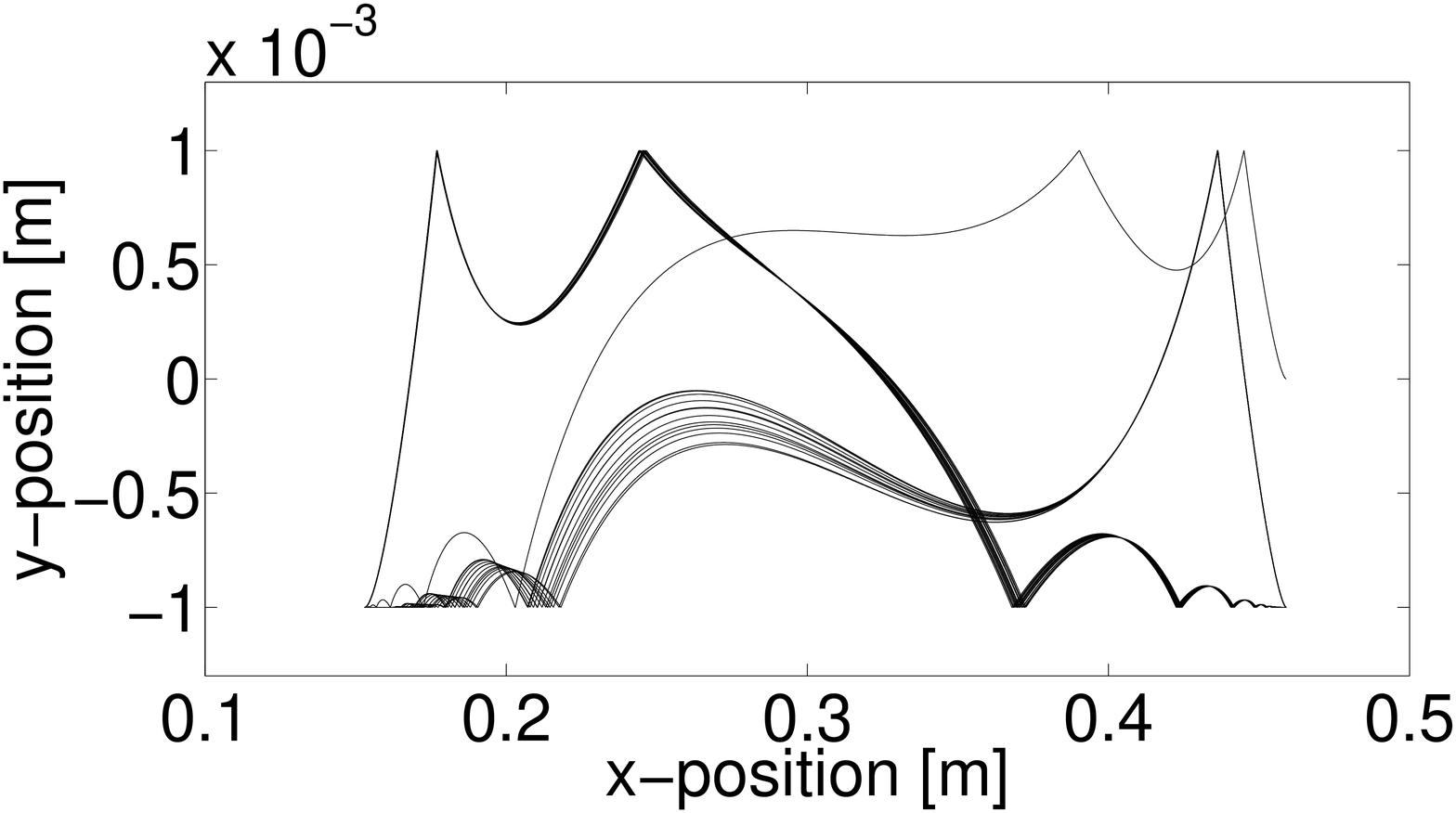}}
  \subfigure[$\epsilon=0.9$]{\includegraphics[width=0.45\columnwidth]{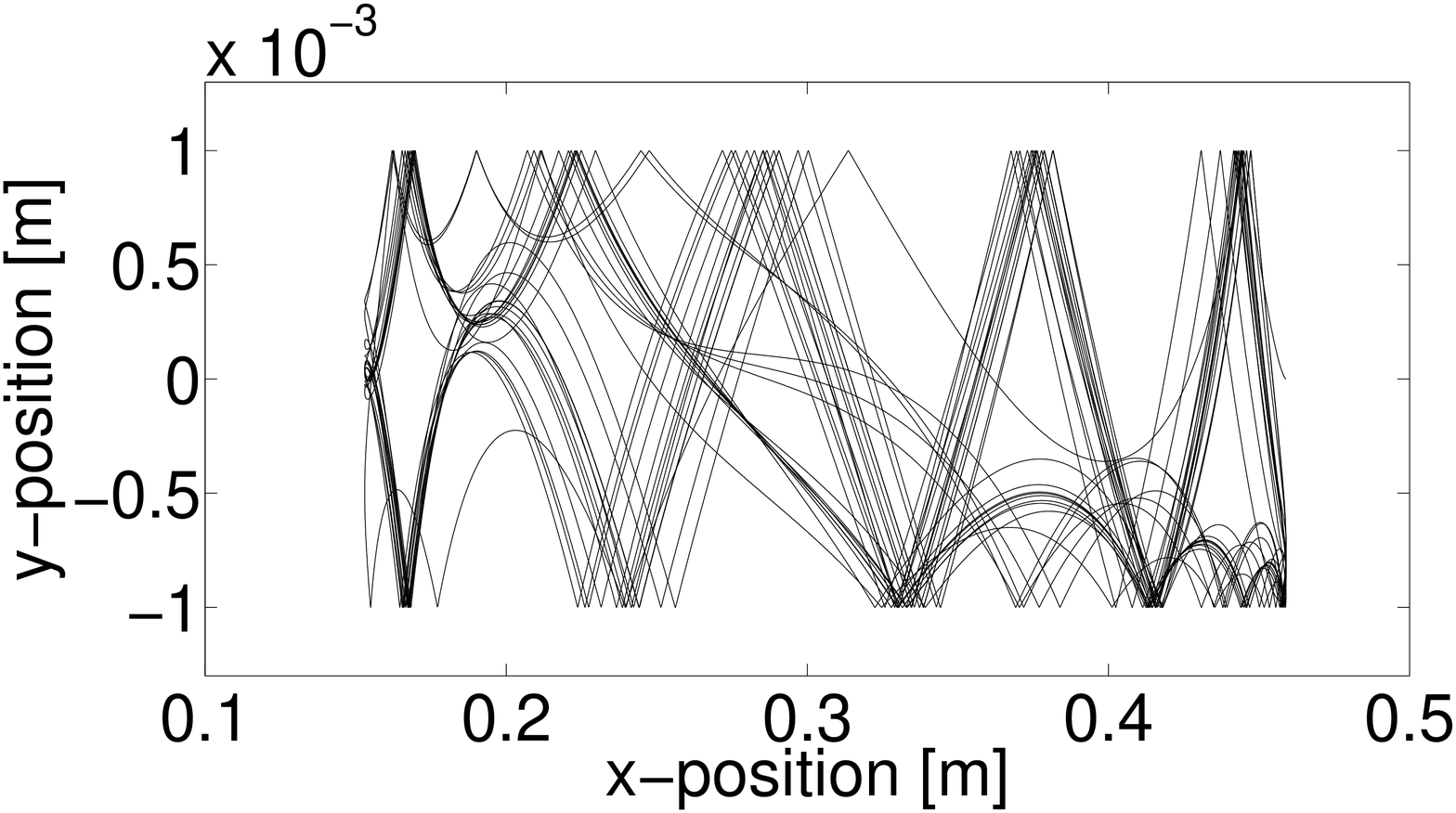}}
  \caption{Movement of the center of gravity of the slider (3) for different coefficients of restitution.}
  \label{fig:resggl}
\end{figure}

In Figure~\ref{fig:resggl}, the results concerning the slider-crank mechanism with unilateral constraints and characteristics as in Table~\ref{tab:slider_crank_characteristics} are presented using a time step size $\Delta t=\unit[10^{-5}]{s}$.
The qualitative behavior is similar to the behavior for Moreau's midpoint rule shown in Figure~\ref{fig:resmoreau}. Especially for high coefficients of restitution, the patterns resemble. The change in the theoretical framework mainly affects persistent contacts, which rarely occur for $\epsilon>0.5$. In contrast for $\epsilon=0.1$, the drift-off effect has a comparatively high influence. The proposed timestepping scheme yields a distinct change in the system's behavior. The drift-off effect does no longer occur and the non-penetration condition is satisfied improving the physical accuracy in comparison to Figure~\ref{fig:resmoreau}.\par
\begin{figure}[hbt]
  \centering
  \includegraphics[width=1.0\columnwidth]{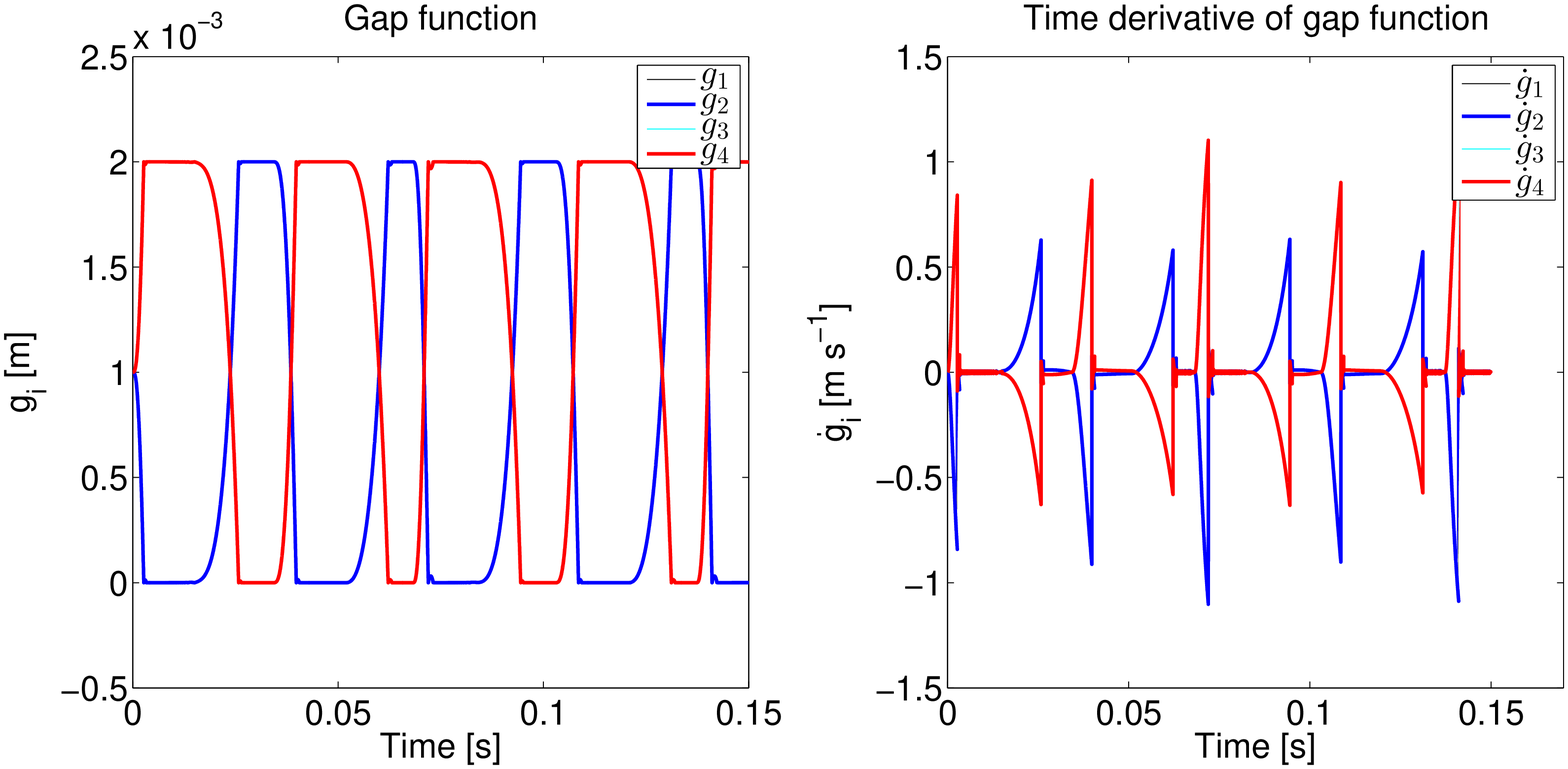}
  \caption{Gap functions and their time derivatives for $\epsilon=0.1$.}
  \label{fig:resgapsggl}
\end{figure}

The development of the gap functions and their time derivatives for $\epsilon=0.1$ is presented in Figure \ref{fig:resgapsggl}.
The drift-off effect has vanished and the gap functions are not negative anymore. The gap velocities are still slightly smaller than zero in time periods where the drift-off effect occurred in the previous simulations. However, this slow trend to increasing permeation is compensated by the second set of Lagrange multipliers enforcing the non-penetration constraint.\par
\begin{figure}[hbt]
  \centering
  \includegraphics[width=1.\textwidth]{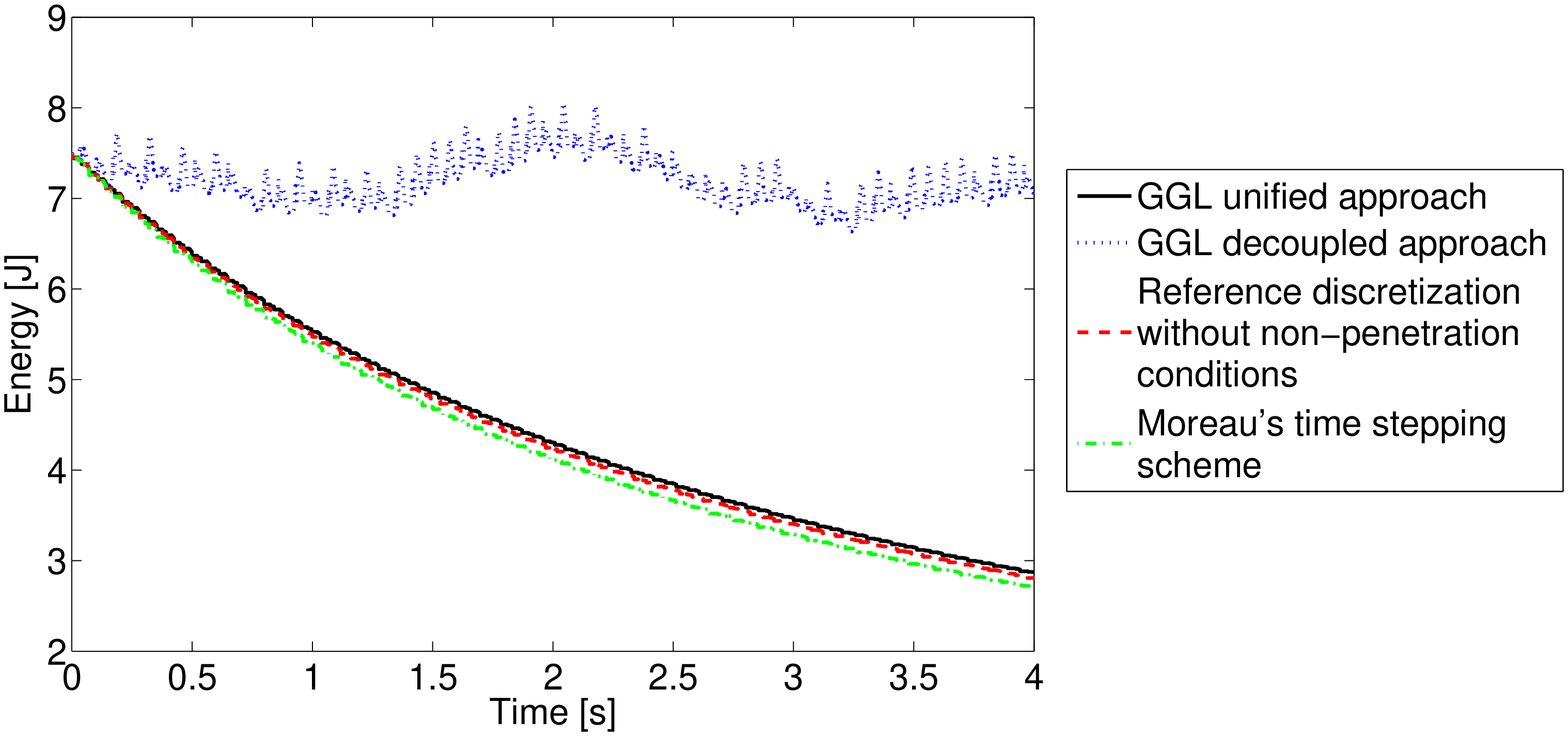}
  \caption{Energy content for different formulations.}
  \label{fig:energycontent}
\end{figure}

To further investigate the physical accuracy of the method, the qualitative development of the entire energy content for Moreau's midpoint rule and the unified Gear-Gupta-Leimkuhler approach is shown in Figure~\ref{fig:energycontent}.
The entire energy content after four seconds of simulation differs slightly. As Moreau's midpoint rule does not ensure the compliance of the constraints, a reference line is shown based on the same algorithm as the unified Gear-Gupta-Leimkuhler approach but only enforcing the impact law like Moreau's midpoint rule and neglecting the non-penetration condition. The reference line is based on \eqref{eq:magnime2} without the last row and without the term $\tilde{\vW}_M^T\left(q_{n+1}\right)\vPsi_{n+1}$ in the first row. The energy development for the reference system is slightly smaller than for the unified Gear-Gupta-Leimkuhler approach and a bit higher than for Moreau's midpoint rule. The proposed approach leads to the same qualitative behavior of the entire energy content and to slightly different quantitative results.\par
Due to the implicit discretization, the computing time increases by a factor of ten in contrast to the explicit Moreau's midpoint rule when using the same time step size. However, we gain a stable discretization which allows comparatively larger time step sizes for stiff problem formulations.  

\section{Conclusion}
Within this work, we propose a timestepping scheme for impacting mechanical systems with unilateral constraints. The new scheme is based on Moreau's midpoint rule and enables to achieve not only compliance of the impact law but also of the non-penetration constraint. It is shown that the decoupled application of the Gear-Gupta-Leimkuhler method for bilateral constraints can be interpreted as a projection to the non-penetration constraint at the end of each time step. As this strategy does not lead to an energy-consistent discretization, our proposition couples position and velocity level with an implicit evaluation of the constraint matrix in the framework of midpoint discretizations. Without significant additional cost, we also approximate the right hand side in the same manner to achieve enhanced stability properties. Adding the active constraints in each time-step by means of the prox~function concept leads to a system of nonsmooth equations which is solved by a Newton scheme.\par
Results from simulations of a slider-crank mechanism with unilateral constraints demonstrate the overcoming of the drift-off effect and the performance of our unified approach. It is reduced concerning the differential index and insofar better conditioned than position level discretizations. Concerning impacting mechanical systems, it is even more important that our scheme is physically consistent due to the impulsive concept and the implicitly incorporated projection. Because of the implicit discretization, the computation time increases significantly in comparison to Moreau's midpoint rule using the same non-controlled time step size. However, for stiff problem formulations, our proposition should result in a more stable discretization and possible larger time step size choices.\par
We do not have applied our scheme to a stiff problem. An analysis concerning numerical issues could be addressed by utilizing backward error analysis~\cite{Hai06}. Using this concept, the interpretation of the induced projection on the non-penetration constraints should be discussed in addition. The relationship to the kinetic metric and appropriate extensions should be studied~\cite{Moel11}. Finally, the consideration of friction taking into account the overdetermined differential algebraic setting according to~\cite{Jay07} would extend our work.

\bibliographystyle{plain}
\bibliography{Literatur}

\begin{thebibliography}{10}

\bibitem{Aca11b}
Vincent Acary and Olivier Bonnefon.
\newblock {Time integration of nonsmooth mechanical systems with unilateral
  contact. Conservation and stability of position and velocity constraints in
  discrete time}.
\newblock In {\em Proceedings of 7th European Nonlinear Oscillation Conference,
  Rome, 24th-29th July 2011}, 2011.

\bibitem{Aca08}
Vincent Acary and Bernard Brogliato.
\newblock {\em {Numerical methods for nonsmooth dynamical systems :
  applications in mechanics and electronics}}, volume~35 of {\em Lecture notes
  in applied and computational mechanics}.
\newblock Springer, Berlin, 1st edition edition, 2008.

\bibitem{Arn09}
Martin Arnold.
\newblock {Numerical methods for simulation in applied dynamics}.
\newblock In Martin Arnold and Werner Schiehlen, editors, {\em Simulation
  Techniques for Applied Dynamics}, number 507 in CISM International Centre for
  Mechanical Sciences, pages 191--246. Springer, Wien, 2009.

\bibitem{Flo10}
Paulo Flores, Remco Leine, and Christoph Glocker.
\newblock Modeling and analysis of planar rigid multibody systems with
  translational clearance joints based on the non-smooth dynamics approach.
\newblock {\em Multibody System Dynamics}, 23:165--190, 2010.

\bibitem{Fueh91}
Claus Führer and Ben Leimkuhler.
\newblock Numerical solution of differential-algebraic equations for
  constrained mechanical motion.
\newblock {\em Numer Math}, 59:55--69, 1991.

\bibitem{Gea85}
Charles~William Gear, Ben Leimkuhler, and G.K. Gupta.
\newblock Automatic integration of {Euler-Lagrange} equations with constraints.
\newblock {\em J Comput Appl Math}, pages 77--90, 1985.

\bibitem{Glo01}
Christoph Glocker.
\newblock {\em {Set-valued force laws in rigid body dynamics : dynamics of
  non-smooth systems}}, volume~1 of {\em Lecture notes in applied and
  computational mechanics}.
\newblock Springer, Berlin, 1st edition edition, 2001.

\bibitem{Hai06}
Ernst Hairer, Christian Lubich, and Gerhard Wanner.
\newblock {\em Geometric numerical integration : structure-preserving
  algorithms for ordinary differential equations}, volume~31 of {\em Springer
  series in computational mathematics}.
\newblock Springer, Berlin, 2nd edition edition, 2006.

\bibitem{Hai10}
Ernst Hairer and Gerhard Wanner.
\newblock {\em Solving Ordinary Differential Equations II: Stiff and
  Differential-Algebraic-Problems}, volume~14 of {\em Springer series in
  computational mathematics}.
\newblock Springer, Berlin, 2nd rev. edition, 1st softcover printing edition,
  2010.

\bibitem{Jay07}
Laurent Jay.
\newblock Specialized partitioned additive {Runge-Kutta} methods for systems of
  overdetermined {DAE}s with holonomic constraints.
\newblock {\em SIAM J Numer Anal}, 45:1814--1842, 2007.

\bibitem{Lei08}
Remco~Ingmar Leine and Nathan van~de Wouw.
\newblock {\em Stability and convergence of mechanical systems with unilateral
  constraints}, volume~36 of {\em Lecture notes in applied and computational
  mechanics}.
\newblock Springer, Berlin, 2008.

\bibitem{Mor99}
Jean~Jacques Moreau.
\newblock Numerical aspects of the sweeping process.
\newblock {\em Comput Methods Appl Mech Engrg}, 177:329--349, 1999.

\bibitem{Moel11}
Michael Möller.
\newblock {\em Consistent integrators for non-smooth dynamical systems}.
\newblock PhD thesis, ETH Zürich, 2011.

\bibitem{Pfe08}
Friedrich Pfeiffer.
\newblock {\em Mechanical system dynamics}, volume~40 of {\em Lecture notes in
  applied and computational mechanics}.
\newblock Springer, Berlin, corr. 2nd printing edition, 2008.

\bibitem{Sch11a}
Thorsten Schindler, Binh Nguyen, and Jeff Trinkle.
\newblock Understanding the difference between prox and complementarity
  formulations for simulation of systems with contact.
\newblock In {\em IEEE/RSJ International Conference on Intelligent Robots and
  Systems, San Francisco, 25th-30th September 2011}, 2011.

\bibitem{Sch12c}
Svenja Schoeder.
\newblock Discussion of the {Gear-Gupta-Leimkuhler} method for unilateral
  contact on the basis of a slider-crank mechanism.
\newblock Semesterarbeit, Technische Universität München, 2012.

\bibitem{Ste11}
David Stewart.
\newblock {\em {Dynamics with inequalities}}.
\newblock SIAM, Philadephia, 2011.

\bibitem{Stu09}
Christian Studer.
\newblock {\em {Numerics of unilateral contacts and friction : modeling and
  numerical time integration in non-smooth dynamics}}, volume~47 of {\em
  Lecture notes in applied and computational mechanics}.
\newblock Springer, Berlin, 2009.

\end{thebibliography}

\end{document}